\documentclass[12pt,a4paper,leqno,twoside]{article}

\usepackage[leqno]{amsmath}
\usepackage{amscd}
\usepackage{amsopn}

\usepackage{amsthm}
\usepackage{amsfonts,bbm}
\usepackage{latexsym}
\usepackage{a4}
\usepackage[all]{xy}
\makeindex
\makeatletter
\@addtoreset{equation}{section}
\makeatother

\newtheorem{lemma}{Lemma}[section]
\newtheorem{proposition}[lemma]{Proposition}

\newtheorem{theorem}[lemma]{Theorem}
\newtheorem{remark}[lemma]{Remark}
\newtheorem{remarks}[lemma]{Remarks}

\newcommand{\real}{\mathbbm{R}}

\newcommand{\nat}{\mathbbm{N}}
\newcommand{\realt}{\mathbbm R^2}
\newcommand{\reald}{\mathbbm R^d}

\renewcommand{\a}{\alpha}
\renewcommand{\b}{\beta}
\newcommand{\ve}{\varepsilon}

\newcommand{\vp}{\varphi}
\newcommand{\g}{\gamma}

\newcommand{\inv}{^{-1}}

\newcommand{\limn}{\lim_{n\to\infty}}

\newcommand{\und}{\quad\mbox{and}\quad}

\begin{document}

\author{Wolfhard Hansen and Nikolai Nikolov%
\thanks{The research of the second author was partially supported 
by a CNRS-grant at the Paul Sabatier University, Toulouse. He is indebted to 
P.J.\ Thomas for stimulating discussions on the subject.}}
\title{One-radius results \\for supermedian functions on
  $\reald$, $d\le 2$ }
\date{}
\maketitle

\begin{abstract}
A classical result states that every lower bounded superharmonic
function on~$\realt$ is constant. In this paper the following
(stronger) one-circle version is proven. If $f\colon \realt\to
(-\infty,\infty]$ is  lower semicontinuous, 
$\liminf_{|x|\to\infty} f(x)/\ln|x|\ge 0$, and, for every $x\in\realt$,
 $1/(2\pi) \int_0^{2\pi}  f(x+r(x)e^{it})\,dt\le f(x)$,  
where $r\colon \realt\to (0,\infty)$ is continuous,
  $\sup_{x\in\realt} (r(x)-|x|)<\infty$, and  
 $\inf_{x\in\realt} (r(x)-|x|)=-\infty$, then $f$ is constant.

Moreover, it is shown that, with respect to   the assumption $r\le
c|\cdot|+M$ on $\real^d$,
there is a striking difference between the restricted
volume mean property for the cases $d=1$ and $d=2$.


\end{abstract}

\section{Introduction and results}

It is a well-known fact that every lower bounded superharmonic
function on~$\real^2$ is constant. We recall that superharmonic
functions~$u$ on~$\realt$ are lower semicontinuous functions on~$\realt$
such that $u>-\infty$, $u\not\equiv \infty$, and, for every circle
$S(x,\rho)$ of center $x\in\realt$ and radius~$\rho>0$, the average
$\sigma_{x,\rho}(u)$ of~$u$
on $S(x,\rho)$ is at most~$u(x)$. In this note, we shall present the
following stronger result (where, as usual, we do not distinguish
between~$\mathbbm C$ and $\realt$).

\begin{theorem} \label{main}
Let $r$ be a strictly positive real function on $\realt$ such that
\begin{itemize}
\item[\rm (i)] $r$ is continuous,
\item[\rm (ii)] $\limsup_{|x|\to\infty}\, (r(x)-|x|)<\infty$,  
\item[\rm (iii)] $\liminf_{|x|\to  \infty}
  (r(x)-|x|)=-\infty$.\footnote{Having (i), properties (ii),(iii)
    are equivalent to $\sup_{x\in\realt}(r(x)-|x|)<\infty$,
    $\inf_{x\in\realt}(r(x)-|x|)=-\infty$, respectively.}
\end{itemize} 
Let $f>-\infty$ be a lower semicontinuous
 numerical function on~$\realt$ such that 
\begin{equation}\label{linf}
\liminf_{|x|\to\infty}  f(x)/\ln|x|\ge 0
\end{equation} 
and $f$  is
\emph{$(\sigma,r)$-supermedian},
that is, 
\begin{equation}\label{sigma}
\sigma_{x,r(x)}(f) :=\frac 1{2\pi}\int_0^{2\pi} f(x+r(x)e^{it})\,dt\
\le\  f(x)\qquad\ \ (x\in\realt).
\end{equation} 
Then $f$ is constant.
\end{theorem}

\begin{remarks} \label{contracting}{\rm
1.  Obviously,  $r\colon \realt\to (0,\infty)$ has the properties
(i) -- (iii), if there exists $L\in (0,1)$ such that, for all $x,y\in\realt$,
\begin{equation}\label{lip}
              |r(x)-r(y)|\le L \,|x-y|.
\end{equation} 
However, assuming only that $|r(x)-r(y)|<
|x-y|$ (which implies (i) and (ii)), even the conclusion breaks down.
Indeed, if $f:=(1-|x|)^+$ and $r:=|x|+2+ (|x|+1)\inv$, then, for all 
$x,y\in \realt$, $\sigma_{x,r(x)}(f)=0\le f(x)$ and the inequality
$|r(x)-r(y)|<|x-y|$ holds, since $(|x|+1)\inv>(|y|+1)\inv$ if
$|x|<|y|$. In fact, none of the properties (i), (ii), (iii) may
be dropped (see Section \ref{counter}). Moreover, since the function
$-\ln(|\cdot|^2+1)$  is superharmonic,   it  is clear that (\ref{linf})
cannot be replaced by $\liminf_{|x|\to\infty}  f(x)/\ln|x|>-\infty$.

2. In Theorem \ref{main},  we may just as well
assume that $f$ does not attain the value $\infty$  
(it  suffices to consider the functions $f_n:=\min\{f,n\}$, $n\in\nat$,
which are $(\sigma,r)$-supermedian provided  $f$ is
$(\sigma,r)$-supermedian; if these functions 
are constant, then $f:=\limn f_n$ is constant).
}
\end{remarks}

Let us assume, for a moment, that $f\colon \realt\to \real$ is
continuous and \emph{$(\sigma,r)$-median}, that is, such that
$$
\sigma_{x,r(x)}(f)=f(x)\qquad (x\in\realt). 
$$
P.C.~Fenton \cite{Fenton79}  
showed that $f$ has to be constant
provided $f$ is lower bounded, $r$ is continuous and, for some
$x_0\in\realt$, the set $\{x\in\real^2\colon r(x)>|x-x_0|\}$ is bounded,
a~requirement which may be replaced by the weaker property (ii)
(see Remark \ref{fenton-kurz}). If $f$ is bounded,
then (ii) alone (without any further assumption on $r$) is sufficient
to conclude that $f$~is~constant (\hbox{\cite[Theorem 1.1]{H-monthly}}, cf.\
also  \cite{H-sigma}).
On the other hand,  there exist  $r\colon \realt\to (0,\infty)$
and a~continuous $(\sigma,r)$-median function~$f$ on $\realt$
such that  $r\le 4(|\cdot|+1)$ and $\min f(\realt)=0$, $\max f(\realt)=1$
 (see \hbox{\cite[Proposition 6.1]{H-sigma}} or \cite[Section 5]{H-monthly}).

An essential step for the strong version \cite[Theorem 1.1]{H-monthly}
of Liouville's theorem   consists in proving that, assuming (ii),  
every lower semicontinuous $(\sigma,r)$-supermedian 
function~$f\ge 0$ on~$\realt$ attains a minimum. It immediately implies that
 constant functions  are the only lower semicontinuous, lower bounded
 functions~$f$  on~$\realt$
which are \emph{$(\lambda,r)$-supermedian}, that is, which have the
property that, for every $x\in \realt$, the average
$\lambda_{x,r(x)}(f)$ of~$f$ on the (closed) disk $B(x,r(x))$ is at
most $f(x)$ (see \cite[Corollary 6.1]{H-monthly}).
We recall that $(\lambda,r)$-supermedian functions are $(\sigma,
\tilde r)$-supermedian for some function $0<\tilde r\le r$
(cf.~\cite[Section~6]{H-monthly}). 
The following result
shows that the existence of a minimum fails, if (ii)  is replaced by   an inequality $r\le c|\cdot|+M$,
where $c>1$.

\begin{proposition} \label{no-min}
Let $c>1$, $M>0$, and $r(x):=\max\{c|x|,M\}$, $x\in\realt$.
Then the functions $r^{-\a}$ on $\realt$
are  $(\sigma,r)$-supermedian {\rm(}and $(\lambda,r)$-supermedian{\rm)}
provided $\a>0$ is~sufficiently small.
\end{proposition}

So, assuming that $r\le c|\cdot|+M$, where $c,M\in (0,\infty)$, a
result of Liouville type for $(\lambda,r)$-supermedian functions 
on~$\realt$ holds if and only if  $c\le 1$. 

On the real line, this will turn out to be strikingly different.
 By the following proposition, such a result of
Liouville type  on $\real$ holds if and
 only if $c\le c_0$, where $c_0\in [2.50, 2.51]$ is~the unique
solution to the equation
$$
(t+1)  \ln(t+1) +  (t-1)  \ln(t-1) =2t
$$
in $(1,\infty)$ (see Section \ref{c0}).

\begin{proposition}\label{d1}
\begin{enumerate} 
\item [\rm 1.] Let $r\colon \real\to (0,\infty)$ and $M>0$ such that 
\begin{equation}\label{critical}
r(x)\le c_0|x|+M \qquad \mbox{ for all } x\in \real\setminus [-M,M].
\end{equation} 
Then every lower semicontinuous $(\lambda,r)$-supermedian
 function $f>-\infty$ on $\real$ which is lower bounded {\rm(}\!or 
satisfies
$\liminf_{|x|\to\infty} f(x)/\ln|x|\ge 0${\rm)}  is constant.
 \item [\rm 2.] If, however, $c> c_0$, $M>0$, and 
 $r:=\max(c|\cdot|,M)$, then
 the function $r^{-\a}$ is $(\lambda,r)$-supermedian
provided $\a>0$ is sufficiently small. 
\end{enumerate} 
\end{proposition}

Finally, let us recall that every bounded harmonic function on
$\reald$, $d\ge 1$, is constant. Hence the results of Liouville type, 
which we discussed until now, are special cases of one-radius results
 for harmonic functions on open sets  $U$. 
The trivial requirement
that $r$ be at most the distance to $U^c$, if $U\ne \reald$, implies
the existence of a real $M>0$ such that
\begin{equation}\label{M}
                  r\le |\cdot|+M
\end{equation} 
(which may justify  considering (\ref{M}) as a natural assumption on $r$
in the case $U=\reald$). 

One radius-results for harmonic functions  have a long history  
(see the survey papers \cite{NV-histmean,H-kouty} and the references therein).
If $U$ is an arbitrary open set  in $\reald$, then every
continuous $(\lambda,r)$-median functions $f\colon U\to \real$ 
admitting a (sub)harmonic
minorant and a (super)harmonic majorant is harmonic (provided
(\ref{M})  holds, if  $U=\reald$). 

Let us now return  to continuous bounded
 $(\sigma,r)$-median functions. If $d=1$, the corresponding result
 fails almost trivially both for~$\real$ and~$(-1,1)$
(in the first case consider $f(x):=\sin x$ and $r(x):=2\pi$,
for the interval see e.g.\ \cite[Section IV.3]{courant-hilbert}, cf.\ also
\cite{HN8}). We already mentioned the positive result for $U=\realt$.
On the unit disk, however, there exists a continuous function 
$0\le f\le 1$ having the one-circle property which is not harmonic
(see \cite{HN-sigma} and \cite{H-expo}). 
The corresponding \emph{general} problems for~$\reald$, $d\ge 3$, are unsolved,
both for the open unit ball and the entire space (if, however, $r$ is Lipschitz
with constant $L\in(0,1)$, the answer is positive \cite[Theorem 2]{Hea}).

\section{Proof of Theorem \ref{main}}

Let  $f\colon \realt\to \real$ be  lower semicontinuous such that 
 (\ref{linf}) holds and $f$ is 
$(\sigma,r)$-supermedian, where $r\colon \realt\to (0,\infty)$
satisfies  (i), (ii), and (iii). 
For all $x\in\realt$ and $\rho>0$, let
$$
        B(x,\rho):=\{y\in\real^2\colon |y-x|\le \rho\} \quad\mbox{ and }\quad
        S(x,\rho):=\{y\in\real^2\colon |y-x|= \rho\}.
$$
By (ii),
there is a real $M>0$ such that
$$
r(x)\le |x|+M \qquad \mbox{ for all } x\in B(0,M)^c.
$$
If $f$ is lower bounded, then, by \cite[Proposition 2.1]{H-sigma}   
(see also \cite{H-monthly}),
\begin{equation}\label{ex-min}
f(x_0)\le f \quad\mbox{ for some }x_0\in B(0,M+2).
\end{equation}  
In fact, a short look at the proof  for (\ref{ex-min})  reveals that it is
valid as well under our weaker assumption (\ref{linf}).
So we may suppose without loss of
generality that $f(0)=0$ and $f\ge 0$ (we can replace~$f$ by
the function $x\mapsto f(x_0+x)-f(x_0)$). Since $f$ is lower semicontinuous, we
know, by (\ref{sigma}), that 
\begin{equation}\label{sigma-con}
      f=0 \mbox{ on }S(x,r(x)), \quad \mbox{ whenever } x\in \real^2
      \mbox{ such that } f(x)=0.
\end{equation} 

Let us  use the technique developed in \cite{Fenton79}.
We define an increasing sequence $(\a_n)$ of continuous real
functions on the unit circle $S:=S(0,1)$ 
by $\a_0:= r(0)$ and
$$
          \a_n(u) := \a_{n-1}(u)+ r(\a_{n-1}(u)u)\qquad (n\in\nat).
$$
By induction, we conclude from (\ref{sigma-con}) that, for all
$u\in S$ and $n\in\nat$,  
\begin{equation}\label{a0}
                 f(\a_n(u)u) = 0.
\end{equation} 
Since $r$ is continuous and strictly positive, we obtain immediately
that $\lim \a_n=\infty$. So~there exists $n\in\nat$ such that
$$
                 \a_n > M    \quad\mbox{ on } S.
$$
For the moment, let us fix $u\in S$. We claim that
\begin{equation}\label{aussen}
         f(\a u)=0 \quad \mbox{ for every } \a \ge \a_n(u).
\end{equation} 
Indeed, suppose that (\ref{aussen}) does not hold. Then there exists 
a maximal real $a$ such that $a\ge \a_n(u)$ and $ f(\a u)=0$
for every  $\a\in [\a_n(u), a]$.
We may join the points $y_0:=a u$ and $y_1:=-\a_n(-u)u$ continuously
by an arc $\gamma:[0,1]\to \real^2$ contained in the set
\begin{equation}\label{Weg}
             \{\a u\colon a\ge \a \ge \a_n(u)\}\cup \{  \a_n (v) v\colon v\in S\}.
\end{equation} 
In particular, $f\circ\gamma=0$ by (\ref{a0}).
Fix $0<\beta\le r(y_0)$ and let $z:= y_0+\beta u$
so that $|z- y_0|\le r(y_0)$. Clearly, $|z-y_1| \ge r(y_1)$,
since the origin is contained in the line segment from $y_1$ 
to $y_0$, $r(y_1)\le |y_1|+ M$, and $a+\beta\ge \a_n(u)>M$.
By continuity of $r$, there exists $s\in[0,1]$ such that 
$|z-\gamma(s)|=r(\gamma(s))$. Since $f(\gamma(s))=0$,
we conclude by (\ref{sigma-con}) that $f(z)=0$. Thus
$a\ge a+r(y_0)$, a contradiction proving (\ref{aussen}).

Fixing $R>0$ such that $\a_n\le R$ on $S$, we therefore  know that
$ f=0$ on $\real^2\setminus B(0,R)$.
Since $\liminf_{|x|\to  \infty} (r(x)-|x|)=-\infty$, there exists
$x\in\real^2$ such that $|x|-r(x)>R$, and hence 
\begin{equation}\label{r-small}
     f=0\quad\mbox{ on }B(x,r(x)).
\end{equation} 
Suppose that there is a point $y\in \real^2$ such that $f(y)>0$. 
We define
$$
     t:=\sup\{s\in [0,1]\colon f(sx+(1-s)y)>0\} 
\quad\mbox{ and }\quad z:=tx+(1-t)y.
$$
Since $r(z)>0$,  there exists a point $\tilde y\in [y,z]$ such that
$|\tilde y- z|<r(z)$ and  $f(\tilde y)>0$. By~(\ref{r-small}),
$|\tilde y- x|>r(x)$. By continuity of $r$, we conclude that there exists 
$\tilde z\in (z,x)$ such that $|\tilde y-\tilde z|=r(\tilde z)$.
So $f(\tilde z)>0$, by (\ref{sigma-con}).
However, by definition of $z$, $f=0$ on $(z,x)$.
Thus  there is no point $y\in \real^2$ such that  
$f(y)>0$, that is, $f$ is identically zero, and the proof of 
Theorem~\ref{sigma} is finished.

\begin{remark}\label{fenton-kurz}{\rm
If $f$ is even continuous and $(\sigma,r)$-median, then
(iii) is not needed to conclude that $f$ is constant.

Indeed, it suffices to observe that (iii) has not been used
to obtain that $f=\inf f(\realt)$ outside a compact set,
and hence $\g:=\sup f(\realt)<\infty$.
Then, just using (i) and (ii), we~get as well that
$\g-f=\inf (\g-f)(\realt)$ outside a compact set,
that is, $f=\sup f(\realt)$ outside a compact set.
Thus $\inf f(\realt)=\sup f(\realt)$, $f$ is constant.
}
\end{remark}
\section{Examples}\label{counter}

Simple examples show that a continuous bounded
$(\sigma,r)$-supermedian function $f$ on $\realt$  may  be non-constant, if
any of the properties (i), (ii), or (iii) of $r$ is violated.

1. Let $f:=(1-|x|)^+$, $x\in\real^2$. 
Taking $r(x):=3$, if $|x|<2$, and $r(x):=1$, if
$|x|\ge 2$, we observe that, of course, property (i), that is, the
continuity of $r$, cannot be omitted (or replaced by lower
semicontinuity). Considering $r:=|x|+2+ (|x|+1)\inv$ we already noted
in Remark \ref{contracting}.1 that (iii) cannot be dropped (of course, for this
purpose, it would be sufficient to take $r(x):=|x|+2$).

2. Finally, let us prove that the conclusion of   Theorem
\ref{main} fails, if property (ii) is omitted.
 For $x\in\real^2$, let
$$
   f(x):=\min\{1,|x_1|\inv\} \quad\mbox{ and }\quad r(x):=6
   \max\{1,x_1^2\}.
$$
Clearly, $0\le f\le 1$,  $f$ and $r$ are continuous functions, 
and $r$ satisfies (iii), since
$r(0,t)=6$ for every $t\in\real$. 
To prove that (\ref{sigma}) holds, we  fix $x\in\real^2$ and define 
$$
a:=\max\{|x_1|,1\}, \qquad         A:=\{y\in\real^2\colon |y_1|\le 2a\}.
$$
Then $f(x)= a\inv$. 
We shall see that
\begin{equation}\label{sig}
   \sigma_{x,r(x)}(A)\le a\inv/2
\end{equation} 
and hence
$$
    \sigma_{x,r(x)}(f)\le
     \sigma_{x,r(x)}(A) +\sup f(\real^2\setminus A)
     \le  a\inv/2+ a\inv/2= f(x).
$$
To prove (\ref{sig}) let $\a$ denote the maximal angle between the $x_2$-axis and the lines
connecting~$x$ with one of the four  points $y\in S(x,r(x))$
satisfying $|y_1|=2a$. Then
$$
     \sigma_{x,r(x)}(A)\le 4 \frac \a {2\pi}=\frac 2\pi \a\le \sin \a\le
     \frac{3a}{r(x)}.
$$
If $|x_1|\le 1$, then $a=1$, $r(x)=6$, and hence
$3a/r(x)=1/2=1/(2a)$.
If $|x_1|>1$, then $a=|x_1|$, $r(x)=6x_1^2$, and hence $3a/r(x)=1/(2|x_1|)=1/(2a)$.
Thus (\ref{sig}) holds

A closer look would reveal that $f$ is $\tilde r$-superharmonic with
respect to a continuous function $0<\tilde r\le r$ satisfying
$  \tilde  r(x)=\tilde r(|x_1|,0)$ and
\begin{equation}\label{optimal}  
 \lim_{t\to \infty} 
\frac{\tilde r(t,0)}{t\ln t}=\frac 2\pi.
\end{equation} 
This follows from the fact that, given $t\ge 1$ and $k\in\nat$, the
point $x:=(t,0)$ and the set $\tilde A:=\{y\in\real^2\colon |y_1|\le kt\}$
satisfy
$$
\int_{\real^2\setminus \tilde A} f \,d\sigma_{x,\rho}
\le \sup f(\real^2\setminus \tilde A)\le \frac 1k f(x) \qquad \mbox{
  for every }\rho>0,
$$
and, for large $\rho$,
$$
\int_{\tilde A} f\,d\sigma_{x,\rho}\sim 4\cdot \frac 1{2\pi\rho}
\int_1^{kt} \frac 1\tau\,d\tau=\frac{2\ln(kt)}{\pi\rho}
\sim \frac 2\pi \cdot\frac{t\ln t}\rho\, f(x)
$$
(which, incidentally, shows that the limit behavior in (\ref{optimal})
is optimal for our function~$f$).

\section{Proof of Proposition \ref{no-min}}

Let $c>1$, $M>0$, and $r:=\max\{c|\cdot|,M\}$ so that, for every $\a>0$,
$$
r^{-\a}(x)=\min\{|cx|^{-\a},M^{-\a}\}, \qquad x\in\realt.
$$
We define
$$
 I(\a):=\frac 1{2\pi}\int_0^{2\pi} |1+ce^{it}|^{-\a}\,dt, \qquad
 \a\ge 0.
$$
Then $I(0)=1$ and
$$
   I'(0)=-\frac 1{2\pi}\int_0^{2\pi}\ln|1+ce^{it}|\,dt=-\ln c<0.
$$
So there exists $\a_0>0$ such that $I<1$ on $(0,\a_0]$.
Let us fix $\a\in (0,\a_0]$ and $x\in \realt$. If~$c|x|>M$, then
$r(x)=c|x|$ and hence
$$
\sigma_{x,r(x)}(r^{-\a})\le \sigma_{x,r(x)}(c^{-\a}|\cdot|^{-\a})=
\frac {c^{-\a}}{2\pi}
\int_0^{2\pi}| x+c|x|e^{it}|^{-\a}\,dt = |cx|^{-\a}I(\a)<r^{-\a}(x).
$$
If $c|x|\le M$, then $r(x)=M$, and hence 
$ \sigma_{x,r(x)}(r^{-\a})\le M^{-\a}=r^{-\a}(x)$. 
Thus $r^{-\a}$ is $(\sigma,r)$-supermedian.

Since $\lambda_{x,\rho}=2\rho^{-2}\int_0^\rho \sigma_{x,s}\, s ds$
(and $\int_0^{2\pi} \ln |1+ce^{is}|\,ds=2\pi\ln 1=0$, if $s\in
(0,1)$), we~obtain similarly that $r^{-\a}$ is
$(\lambda,r)$-supermedian provided $\a>0$ is sufficiently small.

\section{Proof of Proposition \ref{d1}}\label{c0}

Let us define 
$$
\psi(t):=(t+1)\ln(t+1)+(t-1)\ln(t-1)-2t, \qquad 1<t<\infty.
$$
Then $\psi$ is continuous,   $\lim_{t\to 1} \psi(t)=2\ln 2-2<0$,
$\lim_{t\to\infty}\psi(t)=\infty$. Moreover,
$$
  \psi'(t)=\ln(t+1)+\ln(t-1)=\ln(t^2-1),
$$
hence $\psi$ is strictly decreasing on $(0,\sqrt 2\,)$ and strictly increasing
on $(\sqrt 2,\infty)$. So there exists 
$c_0\in (1,\infty)$ such that $\psi(c_0)=0$, 
$$
\psi<0 \mbox{ on } (1,c_0), \quad\mbox{ and }\quad \psi>0 \mbox{ on } (c_0,\infty).
$$
In fact,
$2,50<c_0<2,51$ (since $\psi(2.50)<0$ and $\psi(2.51)>0$).

1.  Let $r\colon \real\to (0,\infty)$ and $M>0$ such that 
$$
r(x)\le c_0|x|+M \qquad \mbox{ for all } x\in \real\setminus [-M,M].
$$
 Let $\vp:=\ln^+(|\cdot|-M)$ (so that $\vp(x)=0$, if
  $-(M+1)\le x\le   M+1$). We claim that there exists $\tilde M>0$
such that, for every $x\in \real\setminus [\tilde M,-\tilde M]$,
\begin{equation}\label{l-mean}
        \lambda_{x,r(x)}(\vp)\le \vp(x).
\end{equation} 
Since $\lim_{x\to \infty} (x/M)\ln\bigl(1-(M/x)\bigr)=-\ln'1=-1$, there
exists $\tilde M\ge 1+c_0+2M$   such that 
\begin{equation}\label{tM}
M\ln (x-M) + c_0x \ln \frac {x-M}x-1>0 \quad\mbox{ for every }x>
\tilde M.
\end{equation} 
For a while, let us fix $ x\in\real\setminus [-\tilde M,\tilde M]$. 
To prove (\ref{l-mean}) we may assume, by symmetry, that $x$ is positive.
Let $y\in (0,c_0x+M)$. If $x-y\ge M+1$,  then 
\begin{equation}\label{y-mean}
  \vp(x-y)+\vp(x+y)\le 2\vp(x),
\end{equation} 
since $\vp$ is concave on $(M+1,\infty)$. If $M+1>x-y\ge -(M+1)$, then 
(\ref{y-mean}) holds, since $\vp(x-y)=0$ and $x+y-M\le x(1+c_0)\le (x-M)^2$.
Therefore (\ref{l-mean}) holds, if~$x-r(x)\ge -(M+1)$.

Let us assume next that $x-r(x)<-(M+1)$. Then
$         t:=(r(x)-M)/x\in (1,c_0]$,
\begin{eqnarray*} 
&&  \int_0^{r(x)\pm x} \vp(s)\,ds=\int_{M+1}^{(t\pm 1)x+M} \vp(s)\,ds\\
&=& \int_1^{(t\pm 1)x}\ln s\,ds=(t\pm 1)x\ln[(t\pm1)x]-(t\pm 1)x +1.
\end{eqnarray*} 
Since $\psi(t)\le \psi(c_0)= 0$,  we hence see  that
$$
\int_{x-r(x)}^{x+r(x)} \vp(s)\,ds=\psi(t)x+2tx\ln x+2\le 2tx\ln x+2,
$$
that is, $r(x)\lambda_{x,r(x)}(\vp)\le tx\ln x+1$.
Thus 
\begin{eqnarray*} 
r(x)(\vp(x)-\lambda_{x,r(x)}(\vp))
&\ge& (tx+M)\ln (x-M)-(tx\ln x+1)\\
&=& M\ln(x-M)+tx\ln \frac {x-M}x-1,
\end{eqnarray*} 
where the right side is positive by (\ref{tM}), since $t\le c_0$. This finishes the proof
of (\ref{l-mean}).

Now let $f>-\infty$ be a lower semicontinuous
$(\lambda,r)$-supermedian function on $\real$ such that 
$\liminf_{|x|\to\infty} f(x)/\ln|x|\ge 0$. To prove that $f$
is constant, we   use ideas from \cite{HN7} and~\cite{H-monthly}.
There exists $x_0\in [-\tilde M,\tilde M]$ such that 
$$
        f(x_0)=\inf f([-\tilde M,\tilde M]).
$$
We intend to show that  $f\ge f(x_0)$ on $\real$. 
Then the lower semicontinuity of $f$ and the inequalities 
$\lambda_{x,r(x)}(f)\le f(x)$, $x\in \real$, will imply that
the set $A:=\{x\in \real\colon f(x)=f(x_0)\}$ is both closed and open, hence $A=\real$, $f=f(x_0)$.

Fixing $\ve>0$, it suffices to prove that
$$
          \tilde f:=f+\ve \vp\ge f(x_0).
$$ 
Obviously, $\tilde f$ is  lower semicontinuous, lower bounded,
and $\lim_{|x|\to \infty} \tilde f(x)=\infty$.
Therefore $\tilde f$ attains a minimum on $\real$, and the non-empty set
$$
         \tilde A:=\{x\in\real\colon \tilde f(x)=\inf \tilde f(\real)\}
$$
is closed. Let $z\in \tilde A$ with minimal absolute value.
If $|z|> \tilde M$, then    $\lambda_{z,r(z)}(\tilde f)\le \tilde f(z)$, and hence
$[z-r(z),z+r(z)]\subset \tilde A$. This is impossible, by our choice
of $z$. Thus
$z\in [-\tilde M,\tilde M]$, and $\tilde f\ge \tilde f(z)\ge f(z)\ge f(x_0)$.

2. Finally, let  $c> c_0$, $M>0$, and 
 $r:=\max(c|\cdot|,M)$. We have to show that
 the function $r^{-\a}$ is $(\lambda,r)$-supermedian
provided $\a>0$ is sufficiently small. 
To that end we define
$$
   \Psi(\b):=(c+1)^\b+(c-1)^\b-2\b c,   \qquad 
                 0<\b<         \infty.
$$
Then $\Psi'(\b)= (c+1)^\b \ln(c+1)+(c-1)^\b\ln(c-1)-2c$. In
particular, $\Psi'(1)=\psi(c)>0$. So there exists $\a\in (0,1)$ such
that $\Psi(1-\a)<0$.

Let us now fix $x\in\real$. If $c|x|\le M$, then $r(x)=M$, 
and hence $ \lambda_{x,r(x)}(r^{-\a})\le M^{-\a}=r^{-\a}(x)$. 
So let us assume that $c|x|<M$ and hence $r(x)=c|x|$. Then
$$
\int_{x-r(x)}^{x+r(x)} |s|^{-\a}\,ds=\frac
1{1-\a}|x|^{1-\a}\bigl((c+1)^{1-\a}+(c-1)^{1-\a}\bigr) \und
2r(x)|x|^{-\a}=2c|x|^{1-\a}.
$$
Since $\Psi(1-\a)<0$, we conclude that $\lambda_{x,r(x)}(|\cdot|^{-\a})\le |x|^{-\a}$
and hence 
$$
\lambda_{x,r(x)}(r^{-\a})\le \lambda_{x,r(x)}(c^{-\a}|\cdot|^{-\a})
\le |cx|^{-\a} =r^{-\a}(x).
$$
Thus $r^{-\a}$ is $(\lambda,r)$-supermedian.

\begin{remark} {\rm
If even $r\le c|\cdot|+M$ for some $c\in(0,c_0)$, then the conclusion
in (1) of~Proposition~\ref{d1} is still valid, if the condition 
$\liminf_{|x|\to \infty} f(x)/\ln |x|\ge 0$ is replaced by the
weaker assumption $\liminf_{|x|\to \infty} f(x)/\ln |x|>-\infty$. Indeed,
by means of the function~$\Psi$, we may then prove that,
for some $\a>0$,  the function $|\cdot|^\a$
(which will replace $\vp$) is~$(\lambda,r)$-supermedian.
}
\end{remark}

{\small \noindent 
Wolfhard Hansen,
Fakult\"at f\"ur Mathematik,
Universit\"at Bielefeld,
33501 Bielefeld, Germany, e-mail:
 hansen$@$math.uni-bielefeld.de}\\
{\small \noindent 
Nikolai Nikolov,
Institute of Mathematics, Acad. G. Bonchev str., block 8, 1113 Sofia,
Bulgaria, e-mail:
nik$@$math.bas.bg}


\begin{thebibliography}{AB1}
\bibliographystyle{alpha}


\bibitem{courant-hilbert}
R.~Courant and D.~Hilbert.
\newblock {\em Methods of mathematical physics. {V}ol. {II}}.
\newblock Wiley Classics Library. John Wiley \& Sons Inc., New York, 1989.
\newblock Partial differential equations, Reprint of the 1962 original, A
  Wiley-Interscience Publication.

\bibitem{Fenton79}
P.~C. Fenton.
\newblock On sufficient conditions for harmonicity.
\newblock {\em Trans. Amer. Math. Soc.}, 253:139--147, 1979.

\bibitem{H-kouty}
W.~Hansen.
\newblock Restricted mean value property and harmonic functions.
\newblock In {\em Potential theory---ICPT 94 (Kouty, 1994)}, pages 67--90. de
  Gruyter, Berlin, 1996.

\bibitem{HN7}
W.~Hansen and N.~Nadirashvili.
\newblock Restricted mean value property on $\mathbbm R^d$, 
$d\leq 2$.
\newblock{\em Expo.  Math.}, 13:93--95, 1995.

\bibitem{H-sigma}
W.~Hansen.
\newblock {A Liouville property for spherical averages in the plane}.
\newblock {\em Math. Ann.}, 319:539--551, 2001.

\bibitem{H-expo}
W.~Hansen.
\newblock Littlewood's one-circle problem, revisited.
\newblock {\em Expo. Math.}, 26(4):365--374, 2008.

\bibitem{H-monthly}
W.~Hansen.
\newblock A strong version of {L}iouville's theorem.
\newblock {\em Amer. Math. Monthly}, 115(7):583--595, 2008.

\bibitem{HN-sigma}
W.~Hansen and N.~Nadirashvili.
\newblock Littlewood's one circle problem.
\newblock {\em J. London Math. Soc. (2)}, 50(2):349--360, 1994.

\bibitem{HN8}
W.~Hansen and N.~Nadirashvili.
\newblock Harmonic functions and averages on shells.
\newblock {\em J. Anal. Math.}, 84:231--241, 2001.

\bibitem{Hea}
D.~Heath.
\newblock Functions possessing restricted mean value properties.
\newblock {\em Proc. Amer. Math. Soc.}, 29:588--595, 1973.

\bibitem{NV-histmean}
I.~Netuka and J.~Vesel{\'y}.
\newblock Mean value property and harmonic functions.
\newblock In {\em Classical and modern potential theory and applications
  ({C}hateau de {B}onas, 1993)}, volume 430 of {\em NATO Adv. Sci. Inst. Ser. C
  Math. Phys. Sci.}, pages 359--398. Kluwer Acad. Publ., Dordrecht, 1994.

\end{thebibliography}
\end{document}